\documentclass{amsart}
\usepackage{amsmath,amssymb,amsthm,mathrsfs}

\newtheorem{thm}{Theorem}[section]
\newtheorem{prop}[thm]{Proposition}
\newtheorem{lem}[thm]{Lemma}
\newtheorem{cor}[thm]{Corollary}

\newcommand{\sF}{{\mathcal F}}

\newcommand{\sW}{{\mathcal W}}

\def\k{\kappa}

\def\t{\tau}
\def\va{\varphi}

\def\om{\Omega}

\def\tht{\Theta}

\def\ts{\times}

\def\iy{\infty}
\def\im{{\rm Im\, }}

\def\kr{{\rm Ker\, }}

\def\ind{{\rm ind\,}}

\def\lg{\langle}
\def\rg{\rangle}
\def\wh{\widehat}
\def\wt{\widetilde}

\def\BC{{\mathbb C}}

\def\BR{{\mathbb R}}

\newcommand{\bpr}{{\noindent\textbf{Proof.}\ \ }}
\newcommand{\epr}{{\hfill $\Box$}}

\newcommand{\nn}{\nonumber}

\newcommand{\ands}{\quad\mbox{and}\quad}

\renewcommand{\theequation}{\arabic{section}.\arabic{equation}}

\begin{document}

\title{The  B\'ezout  equation
 on the right half plane in a Wiener space setting}

\author[G.J. Groenewald]{G.J. Groenewald}
\address{%
G.J. Groenewald, Department of Mathematics, Unit for BMI, North-West University\\
Private Bag X6001-209, Potchefstroom 2520, South Africa}

\email{Gilbert.Groenewald@nwu.ac.za}


\author[S. ter Horst]{S. ter Horst}

\address{%
S. ter Horst, Department of Mathematics, Unit for BMI, North-West University\\
Private Bag X6001-209, Potchefstroom 2520, South Africa}

\email{sanne.terhorst@nwu.ac.za}


\author[M.A. Kaashoek]{M.A. Kaashoek}

\address{%
M.A. Kaashoek, Department of Mathematics,
VU University Amsterdam\\
De Boelelaan 1081a, 1081 HV Amsterdam, The Netherlands}

\email{m.a.kaashoek@vu.nl}

\thanks{The third author gratefully thanks the mathematics department of North-West University, Potchefstroom campus, South Africa, for the   hospitality and  support  during his visit from September 21--October 15, 2015.\\
This work is based on the research supported in part by the National Research Foundation of South Africa (Grant Number 93406).}


\begin{abstract} This paper deals with the  B\'ezout equation $G(s)X(s)=I_m$, $\Re s \geq 0$, in the Wiener space of analytic matrix-valued functions on the right half plane. In particular, $G$ is an $m\ts p$ matrix-valued analytic Wiener function,  where $p\geq m$,  and the solution $X$ is required to be an  analytic Wiener function  of size $p\ts m$.  The set of all solutions is described explicitly in terms of a $p\ts p$ matrix-valued analytic Wiener function $Y$, which has an inverse in the analytic Wiener space, and an associated inner function $\Theta$ defined by $Y$ and the value of $G$ at infinity. Among the solutions, one is identified that minimizes the $H^2$-norm.  A Wiener space version of  Tolokonnikov's lemma plays an important role in the proofs. The results presented  are natural analogs of those obtained for the discrete case in \cite{GtHK1}.
\end{abstract}



\maketitle


\setcounter{section}{0}
\setcounter{equation}{0}
\section{Introduction and main results}\label{sec:intro}
In this paper we deal with the B\'ezout equation $G(s)X(s)=I_m$ on the closed right half plane $\Re s \geq 0$, assuming that  the given function $G$  is of the form
\begin{align}
&G(s)=D+ \int_0^\iy e^{-st}g(t) \,dt \quad( \Re s \geq 0),\nn \\
&\hspace{3cm} \mbox{where}\quad  g\in L^1_{m\ts p}(\BR_+)\cap L^2_{m\ts p}(\BR_+). \label{GWiener}
\end{align}
In particular, $G$ belongs to the analytic Wiener space $\sW_+^{m\ts p}$. We are interested in  solutions $X\in \sW_+^{p\ts m}$, that is,
\begin{equation}\label{solX}
X(s)=D_X+\int_0^\iy e^{-st}x(t) \,dt \quad( \Re s \geq 0), \quad\mbox{where}\quad  x\in L^1_{m\ts p} (\BR_+).
\end{equation}
Throughout $p\geq m$. We refer to the final paragraph of this introduction for a further  explanation of the notation.

With $G$ given by \eqref{GWiener} we associate the Wiener-Hopf  operator $T_G$ mapping $L_p^2(\BR_+)$ into $L_m^2(\BR_+)$   which is defined  by
\begin{equation}\label{WHopfG}
(T_Gh) (t)=D h(t)+\int_0^\iy g(t-\tau)h(\tau) d\tau, \quad t\geq 0\quad (h\in L_m^2 (\BR_+)).
\end{equation}
For  $X$ as in \eqref{solX} we define the Wiener-Hopf  operator $T_X$  mapping $L_m^2(\BR_+)$ into $L_p^2(\BR_+)$  in a similar way, replacing $D$ by $D_X$ and $g$ by $x$. If the B\'ezout equation
\begin{equation}\label{Bez1}
G(s)X(s)=I_m, \quad \Re s \geq 0.
\end{equation}
has a solution $X$ as in \eqref{solX}, then (using the analyticity of $G$ and $X$) the theory of Wiener-Hopf operators (see \cite[Section XII.2]{GGK1} or \cite[Section 9]{BS06}) tells us  that $T_GT_X=T_{GX}=I$, where $I$ stands for the identity operator on $L^2_m(\BR_+)$.  Thus for  the B\'ezout equation \eqref{Bez1} to be solvable  the operator $T_G$ must be surjective or, equivalently, $T_GT_G^*$ must be  strictly positive. We shall see that this condition is also sufficient.

To state our main results,  we  assume that $T_GT_G^*$ is strictly positive.  Then $D = G(\infty)$ is surjective, and hence $DD^*$ is strictly positive too. We introduce two matrices $D^+$ and $E$, of sizes $p\ts m$ and $p\ts (p-m)$, respectively, and a ${p\ts p}$ matrix function $Y$ in $\sW_{+}^{p\ts p}$, as follows:

\begin{itemize}
  \item[\textup{(i)}] $D^+ = D^*(DD^*)^{-1}$, where $D = G(\infty)$;
  \item[\textup{(ii)}] $E$ is an isometry mapping $\mathbb{C}^{p-m}$ into $\mathbb{C}^p$ such that $\im E = \kr D$;
  \item[\textup{(iii)}] $Y$ is the ${p\ts p}$ matrix function given  by
\begin{equation}\label{defY2}
Y(s)=   I_p -  \int_{0}^\iy e^{- s t} y(t) dt,\ \Re s\geq 0, \ \mbox{where $y= T_G^*(T_G T_G^*)^{-1}g$}.
\end{equation}
\end{itemize}
From the definitions  of $D^+$ and $E$ it follows that the $p\ts p$ matrix $\begin{bmatrix} D^+ & E \end{bmatrix}$ is non-singular. In fact
\begin{equation}
\label{idDE}
\begin{bmatrix}D\\ E^*\end{bmatrix}\begin{bmatrix} D^+ & E \end{bmatrix}=
\begin{bmatrix} I_m & 0\\ 0& I_{p-m}\end{bmatrix}.
\end{equation}
As we shall see (Proposition  \ref{prop:HGR} in  Section \ref{sec:TG} below), the fact that the given function $g\in L^1_{m\ts p}(\BR_+)\cap L^2_{m\ts p}(\BR_+)$ implies that a similar result holds true for  $y$. In particular,  $Y\in \sW_+^{p\ts p}$. In what follows $\Xi$ and $\tht$ are the functions defined by
\begin{align}
\Xi(s)&= \Big(I_p -  \int_{0}^\iy e^{- s t} y(t) dt \Big)D^+ = Y(s)D^+,\ \Re s > 0;\label{defXi}\\
\tht(s)&= \Big(I_p -  \int_{0}^\iy e^{- s t} y(t) dt \Big)E = Y(s)E, ,\ \Re s > 0. \label{defTheta}
\end{align}
Since $Y\in \sW_+^{p\ts p}$, we have  $\Xi\in \sW_+^{p\ts m}$ and $\Theta\in \sW_+^{p\ts (p-m)}$. Finally, recall that a function $\om$ in the analytic Wiener space $\sW_+^{k\times r}$ is inner whenever $\om(s)$ is an isometry for each $s \in i \BR$.
We now state our main results.

\begin{thm}\label{mainthm1}
Let $G$ be the ${m\ts p}$ matrix-valued function given by \eqref{GWiener}. Then the equation $G(s) X(s)  = I_m$, $\Re s > 0$,  has a solution $X\in \sW_+^{p\ts m}$  if and only if $T_G$ is right  invertible. In  that case  the function $\Xi$ defined by \eqref{defXi} is a particular  solution and the set of all  solutions $X\in \sW_+^{p\ts m}$ is given by
\begin{equation}\label{sol8}
X(s)= \Xi(s) +  \tht(s)Z(s) ,\quad \Re s > 0,
\end{equation}
where $\Xi$ and $\tht$ are  defined by \eqref{defXi} and \eqref{defTheta}, respectively, and the free parameter $Z$ is an arbitrary  function in $\sW_+^{(p-m)\times m}$.  Moreover, the function $\tht$ belongs to $\sW_+^{p\ts (p-m)}$ and is inner. Furthermore,  the solution $\Xi$ is the minimal $H^2$ solution in the following sense
\begin{align}
& \|X(\cdot) u\|_{H_p^2}^2=\|\Xi(\cdot)u\|_{H_p^2}^2+\|Z(\cdot)u\|_{H_{p-m}^2}^2, \nonumber \\
&\hspace{3cm} \mbox{where   $u\in \BC^m$ and $Z\in \wt{\sW}_+^{\,(p-m)\times m}$}. \label{mincond}
\end{align}
\end{thm}

In the above theorem, for any positive integer $k$, $H_k^2=H_k^2(i\BR)$ is the Hardy space of $\BC^k$-valued functions on the right half plane given by $H_k^2(i\BR)=JL_k^2(\BR_+)$, where $J$ is the unitary operator  defined by
\begin{equation}
\label{defJ}
J=\frac{1}{\sqrt{2\pi}}\sF:  L_k^2(\BR)\to  L_k^2(i\BR)
\end{equation}
with $\sF$ being  the Fourier transform mapping $L_k^2(\BR)$ onto $L_k^2(i\BR)$. Moreover, $Z\in \wt{\sW}_+^{\,(p-m)\times m}$   means that
\begin{align*}
Z(s)&=D_Z+\int_0^\iy e^{-st} z(t)\,dt \quad  (\Re s\geq 0),\quad  \mbox{where}\\
&\hspace{1cm}\mbox{$D_Z$ is a $(p-m)\ts p$ matrix and } z\in L_{(p-m)\ts p}^1(\BR_+) \cap  L_{(p-m)\ts p}^2(\BR_+).
\end{align*}
See the final part of this introduction for further information  about  the used notation, in particular, see \eqref{four} for the definition of the Fourier transform $\sF$.

The second  theorem is a variant of the Tolokonnikov  lemma  \cite{Tol81} in  the present setting. The result emphasizes the central role of the function $Y$.

\begin{thm}\label{mainthm2}
Assume  $T_G T_G^*$ is strictly positive, and let  $Y$ be the matrix function defined  by \eqref{defY2}. Then $Y$  belongs to the Wiener space $\sW_+^{p\ts p}$, $\det Y(s)\not = 0$  whenever   $\Re s \geq  0$, and  hence  $Y$ is invertible in $\sW_+^{p\ts p}$.  Furthermore, the $p\ts p$ matrix function
\begin{equation}\label{out44}
 \begin{bmatrix}
  G(s) \\
  E^* Y(s)^{-1}  \\
\end{bmatrix}, \quad \Re s\geq 0,
\end{equation}
is   invertible  in the Wiener algebra $\sW_+^{p\ts p}$ and
its inverse is given by
\begin{equation}\label{out4}
 \begin{bmatrix}
  G(s) \\
  E^* Y(s)^{-1}  \\
\end{bmatrix}^{-1}=
Y(s) \begin{bmatrix} D^+ &  E   \end{bmatrix}
=\begin{bmatrix} \Xi(s)&  \tht(s) \end{bmatrix}, \quad \Re s\geq 0.
\end{equation}
\end{thm}

\medskip
The literature on the B\'ezout equation and the related corona problem is extensive, starting with   Carleson's corona theorem \cite{Carl62} (for the case when $m=1$) and Fuhrmann's extension to the matrix-valued case \cite{Fuhr68}, both in a   $H^\iy$ setting.  The topic  has beautiful connections with operator theory (see  the books \cite{Helton87},   \cite{Nikol86},  \cite{Nikol02},  \cite{Peller03},  and the more recent papers \cite{Treil04}, \cite{TW05}, \cite{TrentZh06}).  Rational matrix equations of the form \eqref{Bez1}  play an important role in solving systems and control theory problems, in particularly, in problems involving  coprime factorization, see, e.g., \cite[Section 4.1]{Vidya85},  \cite[Section A.2]{GL95}, \cite[Chapter 21]{Zhou96}. For more recent work see  \cite{GuB04} and \cite{GuL03}, and  \cite[page 3]{MS08} where it is proved that  the scalar analytic Wiener algebra is a pre-B\'ezout  ring.  For matrix  polynomials, the equation \eqref{Bez1} is closely related to the Sylvester  resultant; see, e.g., Section 3 in \cite{GKLe08} and the references in that paper.

The present paper is inspired by  \cite{FKR4} and \cite{GtHK1}.  The  paper \cite{FKR4} deals with equation \eqref{Bez1} assuming the matrix function $G$ to be  a stable rational matrix function, and the solutions are required to be  stable rational matrix functions as well. The comment in the final paragraph of  \cite[Section 2]{FKR4} was the starting point for our analysis.  The paper \cite{GtHK1} deals with the discrete case (when the right half plane  is replaced by the open unit disc). Theorems \ref{mainthm1} and \ref{mainthm2} are the continuous analogue of Theorem 1.1 in \cite{GtHK1}.  The absence of an explicit formula for the function $Y^{-1}$ in the present setting makes the proofs  more complicated than those in \cite{GtHK1}.

 The paper consists of  five sections,  including the present introduction and an appendix. Section \ref{sec:TG}, which deals with the right invertibility  of the operator $T_G$, has an auxiliary  character. Theorem \ref{mainthm2} is proved in   Section \ref{sec:prthm2},  and  Theorem~\ref{mainthm1} in Section \ref{sec:prthm1}. The Appendix, Section \ref{sec:App}, contains a number of auxiliary results involving the Lebesgue space $L^1(\mathbb R) \cap L^2(\mathbb R)$ and its vector-valued counterpart, which are collected together simply for the convenience of the reader and contains no significantly new material.

\medskip
\paragraph{Notation and terminology.} We conclude this section with some notation and terminology. Throughout, a linear  map $A:\BC^r \to \BC^k$ is identified with the $k\ts r$ matrix of $A$ relative to the standard orthonormal bases in $\BC^r$ and $\BC^k$. The space  of all ${k\ts r}$ matrices with entries  in $L^1(\mathbb{R})$ will be denoted by  $L_{k\ts r}^1(\BR)$.  As usual $\widehat{f} $ denotes the Fourier transform of $f \in L_{k\ts r}^1(\BR)$, that is,
\begin{equation}\label{four}
\widehat{f}(s) = (\sF f)(s) = \int_{-\infty}^\infty e^{- s t} f(t) dt, \quad s\in i\BR.
\end{equation}
Note that $\wh{f}$ is continuous on the extended imaginary axis $i\BR\cup \{\pm\, i \iy\}$, and  is zero at $\pm\, i\iy$  by  the Riemann-Lebesgue lemma.  By $\sW^{k\ts r}$ we denote the  \emph{Wiener space} consisting of   all  $ k\ts r$ matrix functions  $F$ on the imaginary axis of the form
\begin{align}
&F(s)=D_F+\wh{f}(s), \ s\in i\BR, \mbox{where $f \in L_{k\ts r}^1(\mathbb{R})$ and}\nonumber\\
&\hspace{6cm}\mbox{$D_F$ is a constant matrix.}\label{defF}
\end{align}
Since $\wh{f}$ is continuous on the extended imaginary axis and  is zero at $\pm i\iy$, the function $F$ given by \eqref{defF} is  also continuous on the extended imaginary axis  and the  constant matrix $D_F$  is  equal  to the value of $F$ at infinity. We write $\sW_+^{k\ts r}$ for the space of all $F$ of the form \eqref{defF} with the additional property that  $f$ has its support in $\BR_+=[0, \iy)$, that is, $f$ is equal to zero on $(-\iy, 0)$. Any function $F\in\sW_+^{k\ts r}$ is analytic and bounded on the open right half plane. Thus any  $F\in \sW_+^{k\ts r}$ is  a matrix-valued  $H^\iy$ function. Finally, by $\sW_{-,0}^{k\ts r}$ we denote the Wiener space consisting of   all $F$ of the form \eqref{defF} with the additional property that $D_F=0$ and  $f$ has its support in $(-\iy, 0]$.  Thus we have the following direct sum decomposition:
\begin{equation}
\label{Wdecom}
\sW^{k\ts r}=\sW_+^{k\ts r} \dot{+}\sW_{-,0}^{k\ts r}.
\end{equation}
We write $F\in \wt{\sW}_+^{\,k\ts r}$ if  the function $f$ in  \eqref{defF} belongs to $L_{k\ts r}^1(\BR_+)\cap L_{k\ts r}^2(\BR_+)$. Similarly, $F\in \wt{\sW}_{-,0}^{\,k\ts r}$ if $f\in L_{k\ts r}^1(\mathbb{R_-})\cap L_{k\ts r}^2(\mathbb{R_-})$  and $D_F=0$.

Let  $F\in \sW^{k\ts r}$ be given by \eqref{defF}. With $F$ we associate  the Wiener-Hopf operator $T_F$ mapping  $L^2_r(\BR_+)$ into $L^2_k(\BR_+)$. This operator (see \cite[Section XII.2]{GGK1})  is defined by
\begin{equation}\label{WHopfF}
(T_Fh) (t)=D_F h(t)+\int_0^\iy f(t-\tau)h(\tau) d\tau, \quad t\geq 0\quad (h\in L^2_r(\BR_+)).
\end{equation}
The orthogonal complement of  $H_k^2(i\BR)=JL^2_k(\BR_+)$, with $J$ as in \eqref{defJ},  in  $L_k^2(i\BR)$ will be denoted by  $K_k^2(i\BR)$. If $F\in \wt{\sW}_+^{\,k\ts r}$, then for each $u\in \BC^r$ the function $F(\cdot)u$ belongs to $H_k^2 (i \BR)$. Similarly,  $F(\cdot)u$ belongs to $K_k^2 (i \BR)$  if  $F\in \wt{\sW}_{-,0}^{\,k\ts r}$.

Finally, for $f\in L_{k\ts r}^1(\BR)$ and  $g\in L_{r\ts m}^1(\BR)$ the \emph{convolution product} $f\star g$ is the function in $L_{k\ts m}^1(\BR)$, see \cite[Section 7.13]{Ru66}, given by
\begin{equation}\label{conprod}
 (f\star g) (t)=\int_{-\iy}^\iy f(t-\t)g(\t)\, d\t \hspace{.2cm} \mbox{  a.e.  on $\BR$}.
\end{equation}

\setcounter{equation}{0}
\section{Right invertibility of $T_G$}\label{sec:TG}
In this section $G\in \sW_+^{m\ts p}$, where $G$  is given by \eqref{GWiener} and $p\geq m$.  We already know that the B\'ezout equation \eqref{Bez1}  having  a solution $X$ in $\sW_+^{p\ts m}$ implies that $T_G$ is right invertible or, equivalently, $T_GT_G^*$ is strictly positive; see the paragraph containing  formula \eqref{Bez1}.

In this section we present  an auxiliary result that will be used to prove our  main theorems. For this purpose we need the  $m\ts m$ matrix-valued function $R$ on the imaginary axis defined  by $R(s) = G(s) G({s})^*$, $s\in i\BR$. It follows that $R\in \sW^{m\ts m}$. By  $T_R$ we denote  the corresponding  Wiener-Hopf operator acting on  $L^2_m(\BR_+)$.  Thus
\begin{align*}
(T_Rf)(t)&=DD^*f(t)+ \int_0^\iy  r(t-\t)f(\t)\, d\t,  \quad 0\leq t<\iy, \\
  &  \mbox{with}\  r(t)=Dg^*(t)+g(t)D^*+  \int_{-\iy}^\iy g(t-\t)g^*(\t)\,d\t, \quad t\in \BR.
\end{align*}
Here $g^*(t)=g(-t)^*$ for $t\in \BR$.  It is well-known (see, e.g., formula (24) in Section XII.2 of  \cite{GGK1}) that
\begin{equation}\label{thank}
T_R =   T_GT_G^* +  H_GH_G^*.
\end{equation}\label{HankelG}
Here $H_G$ is the Hankel operator mapping $L_p^2(\BR_+)$ into $L_m^2(\BR_+)$ defined by $G$,  that is,
\begin{equation}\label{defHG}
(H_G f)(t)=\int_0^\infty g(t+\tau)f(\tau) d\tau, \quad f\in L^2_p(\BR_+).
\end{equation}

We shall prove the following proposition. For the case when $G$ is a rational matrix function,  the first part (of the ``if and only if'' part) of the proposition is covered by Lemma 2.3 in \cite{FKR4}.  The proof given in \cite{FKR4} can also be used in the present setting.  For the sake of completeness we  include a proof of the first part.

\begin{prop}\label{prop:HGR} Let $G$ be given by \eqref{GWiener}. Then the operator $T_G$ is right  invertible if and only if     $T_R$  and   $I-H_G^*T_R^{-1}H_G$ are both invertible operators. In that case the  inverse of $T_GT_G^*$ is given by
\begin{equation}\label{TT*inv}
(T_GT_G^*)^{-1}  =T_R^{-1}+T_R^{-1}H_G(I-H_G^*T_R^{-1}H_G)^{-1}H_G^*T_R^{-1}.
\end{equation}
Furthermore,
\begin{itemize}
\item[\textup{(a)}] $(T_GT_G^*)^{-1}$  maps  $L_m^1(\BR_+)\cap L_m^2(\BR_+)$ in a one-to-one way onto itself;
\item[\textup{(b)}] the function $y$ defined by $y=T_G^*(T_G T_G^* )^{-1}g$ belongs to $L^1_{p\ts p}(\BR_+)\cap L^2_{p\ts p}(\BR_+)$, in particular, the function $Y$ given by \eqref{defY2} is in $\wt{\sW}_+^{p\ts p}$.

\end{itemize}
\end{prop}
\bpr  We split the proof into four  parts. In the first  part we assume that $T_G$ is right invertible, and we show that  $T_R$  and   $I-H_G^*T_R^{-1}H_G$ are both invertible operators and that the  inverse of $T_GT_G^*$ is given by \eqref{TT*inv}.  The second part  deals with the reverse implication. Items (a) and (b) are proved in the last two parts.

\smallskip
\noindent\textsc{Part 1.}  Assume    $T_G$ is  right invertible. Then  the operator   $T_G T_G^*$ is strictly positive. According to \eqref{thank} we have $T_R\geq T_G T_G^*$, and hence  $T_R$ is also strictly positive. In particular, $T_R$ is invertible.   Rewriting \eqref{thank}   as   $T_G T_G^* = T_R - H_GH_G^*$, and multiplying  the latter identity from the left and from the right by $T_R^{-1/2}$  shows that
\begin{equation}
\label{nonneg}
T_R^{-1/2}T_G T_G^* T_R^{-1/2}=I-T_R^{-1/2}H_GH_G^*T_R^{-1/2}.
\end{equation}
Hence $I-T_R^{-1/2}H_GH_G^*T_R^{-1/2}$ is strictly positive which shows that $H_G^*T_R^{-1/2}$ is a strict contraction. But then $H_G^*T_R^{-1}H_G=\big(H_G^*T_R^{-1/2}\big)\big(H_G^*T_R^{-1/2}\big)^*$ is also a strict contraction, and thus the operator $I- H_G^*T_R^{-1}H_G$ is strictly positive. In particular, $I- H_G^*T_R^{-1}H_G$  is invertible.  Finally, since  $T_G T_G^* = T_R - H_GH_G^*$, a usual Schur complement type of argument (see, e.g.,   Section 2.2 in \cite{BGKR08}),  including  the  well-known inversion formula
\[
(A - B C)^{-1} = A^{-1} + A^{-1}B(I - C A^{-1} B)^{-1} C A^{-1},
\]
then shows that $(T_G T_G^*)^{-1}$ is given by \eqref{TT*inv}.

\smallskip
\noindent\textsc{Part 2.} In this part we assume that $T_R$  and   $I-H_G^*T_R^{-1}H_G$ are both invertible operators, and we show that $T_G$ is right invertible. According to \eqref{thank}  the operator $T_R$ is  positive. Since we assume  $T_R$ to be invertible, we conclude that $T_R$ is strictly positive.  Rewriting \eqref{thank}   as   $T_G T_G^* = T_R - H_GH_G^*$, and multiplying  the latter identity from the left and from the right by $T_R^{-1/2}$  we obtain the identity \eqref{nonneg}. Hence $I-T_R^{-1/2}H_GH_G^*T_R^{-1/2}$ is positive which shows that $H_G^*T_R^{-1/2}$ is a contraction. But then $H_G^*T_R^{-1}H_G=\big(H_G^*T_R^{-1/2}\big)\big(H_G^*T_R^{-1/2}\big)^*$ is also a contraction, and thus the operator $I- H_G^*T_R^{-1}H_G$ is positive. By assumption  $I- H_G^*T_R^{-1}H_G$ is  invertible. It follows that $I- H_G^*T_R^{-1}H_G$ is strictly positive, and hence  $T_R^{-1/2}H_G$ is a strict contraction. But then the same   holds true for $T_R^{-1/2}H_GH_G^*T_R^{-1/2}$.  This implies that $I-T_R^{-1/2}H_GH_G^*T_R^{-1/2}$ is strictly positive, and \eqref{nonneg} shows that  $T_G$ is right invertible.

\smallskip
\noindent\textsc{Part 3.} In this part we prove item (a). Observe that $g\in L_{m\ts p}^1(\BR)$ implies that  $T_G$ maps $L_p^1(\BR_+)\cap L_p^2(\BR_+)$ into
$L_m^1(\BR_+)\cap L_m^2(\BR_+)$. Since $g^*\in L_{p\ts m}^1(\BR)$ and
\[
(T_G^*f)(t)= D^*f(t)+\int_{0}^\iy  g^*(t-\t)f(\t)\,d\t, \quad 0\leq t<\iy,
\]
the operator   $T_G^*$  maps $L_m^1(\BR_+)\cap L_m^2(\BR_+)$ into $L_p^1(\BR_+)\cap L_p^2(\BR_+)$. Thus $T_GT_G^*$ maps $L_m^1(\BR_+)\cap L_m^2(\BR_+)$ into itself. We have to show that the same holds true for its inverse. To do this we apply    Lemmas \ref{lemR} and \ref{lem4}.

Lemma \ref{lemR} tells us that  $T_R^{-1}$ maps  $L_m^1(\BR_+)\cap L_m^2(\BR_+)$ in a one-to-one way onto itself.  This allows us to apply    Lemma \ref{lem4} with
\[
Q=T_R^{-1}, \quad H=H_G\ands \tilde{H}=H_G^*.
\]
Recall that $H_G$ is a Hankel operator, see \eqref{defHG},  and  $H_G^*$ is also a Hankel operator, in fact
\[
(H_G^*f)(t)=\int_{0}^\iy  g(t+\t)^*f(\t)\,d\t, \quad 0\leq t<\iy.
\]
Since $I- H_G^*T_R^{-1}H_G$ is invertible, Lemma \ref{lem4} then shows that   $I-H_G^*T_R^{-1}H_G$ maps $L_p^1(\BR_+)\cap L_p^2(\BR_+)$ in a one-to-one way onto itself, and hence the same holds true for  its inverse $(I-H_G^*T_R^{-1}H_G)^{-1}$.  To complete the proof  of item (a) note that $H_G^*T_R^{-1}$ maps $L_m^1(\BR_+)\cap L_m^2(\BR_+)$ into $L_p^1(\BR_+)\cap L_p^2(\BR_+)$, and $T_R^{-1}H_G$ maps $L_p^1(\BR_+)\cap L_p^2(\BR_+)$ into $L_m^1(\BR_+)\cap L_m^2(\BR_+)$. But then \eqref{TT*inv}  shows that $(T_GT_G^*)^{-1}$  maps $ L^1_{m}(\BR_+)\cap L^2_{m}(\BR_+)$ into itself. To see that $(T_GT_G^*)^{-1}$ is one-to-one on $L^1_{m}(\BR_+)\cap L^2_{m}(\BR_+)$ and maps $L^1_{m}(\BR_+)\cap L^2_{m}(\BR_+)$ onto itself, one can follow the same argumentation as in the last part of the proof of Lemma \ref{lemR}.

\smallskip
\noindent\textsc{Part 4.} In this part we prove item (b).   Since $g$ belongs to $L^1_{m\ts p}(\BR_+)\cap L^2_{m\ts p}(\BR_+)$,   item (a) tells us that  $f:=(T_GT_G^*)^{-1}g$ also belongs to $L_{m\ts p}^1(\BR_+)\cap L_{m\ts p}^2(\BR_+)$.  We already have  seen (in the first paragraph of the  previous part) that  $T_G^*$  maps $L_m^1(\BR_+)\cap L_m^2(\BR_+)$ into $L_p^1(\BR_+)\cap L_p^2(\BR_+)$.  It follows that $y=T_G^*f$ belongs to $L^1_{p\ts p}(\BR_+)\cap L^2_{p\ts p}(\BR_+)$, as desired. \epr

\setcounter{equation}{0}
\section{The functions $Y$  and $\tht$, and proof of Theorem \ref{mainthm2}}\label{sec:prthm2}

We begin with three lemmas involving the functions $Y$ and $\tht$ defined by \eqref{defY2} and  \eqref{defTheta}, respectively. From Proposition \ref{prop:HGR}, item (b),   and  \eqref{defTheta} we know that $Y\in \sW_+^{p\ts p}$ and $\tht\in  \sW_+^{p\ts (p-m)}$; see also the paragraph preceding Theorem \ref{mainthm1}.

\begin{lem}\label{lem:GY} Assume that  $T_G$ is right invertible, and  let $Y\in \sW_+^{p\ts p}$  be the function defined by \eqref{defY2}. Then
\begin{equation} \label{propY1}
 G(s)Y(s) = D, \qquad \Re s > 0.
\end{equation}
 \end{lem}
 \bpr  To prove \eqref{propY1} note that $T_Gy=T_GT_G^*(T_G T_G^*)^{-1}g=g$. Since the functions $g$ and $y$ both have their support in $\BR_+$, the   identity $T_Gy=g$ can be rewritten as  $Dy+ g\star y=g$, where  $\star$ is the convolution product of matrix-valued functions with entries in $L^1(\BR)$; see \eqref{conprod}. Thus
\begin{equation}
\label{convol1}
Dy(t)+ (g\star y)(t)=Dy(t)+\int_{-\iy}^\iy g(t-\t)y(\t) dt=g(t), \quad t\in \BR.
\end{equation}
Next use that  the Fourier transform of a convolution product is just the product of the Fourier transforms of the functions in the convolution product. Thus taking Fourier transforms in \eqref{convol1} yields  $D\wh{y}+\wh{g}\wh{y}=\wh{g}$. The latter identity can be rewritten as $G\wh{y}=\wh{g}$.  Hence, using the definition of $Y$ in \eqref{defY2}, we obtain
\[
G(s)Y(s)=G(s)\Big(I_p-\wh{y}(s)\Big)=G(s)-\wh{g}(s)=D.
\]
This proves  \eqref{propY1}.\epr

\begin{lem}\label{lem:Theta}
Assume that $T_G$ is right invertible. Then the function $\tht$ defined by \eqref{defTheta} belongs to $\sW_+^{p\ts (p-m)}$ and  is an inner function, that is,  $\tht(s)$ is an isometry for each $s\in i\BR$ and at infinity.
\end{lem}

\bpr
We already know that $\Theta\in\sW_+^{p\ts (p-m)}$. To prove  that $\tht$ is inner, let $y= T_G^*(T_G T_G^*)^{-1}g$ as in \eqref{defY2}, and put $f=(T_GT_G^*)^{-1}g$.  Thus $f\in L^1_{m\ts m}(\BR_+)$,  by Proposition \ref{prop:HGR} (b), and $y=T_G^* f $.  The latter can be rewritten as
\[
y(t)=D^*f(t)+\int_0^\iy g^*(t-\t)f(\t)\, d\t, \quad t\geq 0.
\]
Note that   $g^*(t)=g(-t)^*$,  and hence $g^*$ has its support in $(-\iy, 0]$. Therefore
\begin{align}
y(t)&=D^*f(t)+\int_{-\iy}^\iy g^*(t-\t)f(\t)\, d\t, \quad t\in \BR.\label{yf4}
\end{align}
Put
\begin{equation}\label{defrho}
\rho(t)=\left\{
\begin{array}{cl}
0& \mbox{when $t\geq 0$},\\
(g^*\star f)(t) & \mbox{when $t< 0$.}
\end{array}\right.
\end{equation}
Using the definition of the convolution product $\star$, see    \eqref{conprod}, we  can  rewrite \eqref{yf4} as
\begin{align*}
y(t)&=D^*f(t)+(g^*\star f)(t)-\rho(t)\quad t\in \BR.
\end{align*}
Taking Fourier transforms we  obtain
\[
\wh{y}(s)=D^* \wh{f}(s)+\wh{g^*}(s)\wh{f}(s)- \wh{\rho}(s)
=G(s)^*  \wh{f}(s)-\wh{\rho}(s), \quad s\in i\BR.
\]
Hence, we have
\begin{align}
Y(s)&=I-G(s)^*  \wh{f}(s) + \wh{\rho}(s), \quad s\in i\BR. \label{eqY3}
\end{align}
Now let us compute  $\tht(s)^*\tht(s)=E^*Y(s)^*Y(s)E$  for $s\in i\BR$. We have
\begin{align*}
&E^*Y(s)^*Y(s)E=\\
&\hspace{1cm} =E^*Y(s)^*E-E^*Y(s)^*G(s)^*  \wh{f}(s)E+E^*Y(s)^* \wh{\rho}(s)E\\
&\hspace{1cm} =E^*Y(s)^*E+E^*Y(s)^* \wh{\rho}(s)E\\
&\hspace{3.4cm} \mbox{(because $G(s)Y(s)E=0$ by \eqref{propY1} and $DE=0$)}\\
&\hspace{1cm} =E^*E-E^*\wh{y}(s)^*E+E^*Y(s)^* \wh{\rho}(s)E\\
&\hspace{1cm} =I_{p-m}- \om(s).\end{align*}
Here $\om(s)= E^*\wh{y}(s)^*E-E^*Y(s)^* \wh{\rho}(s)E$. Note that the functions $\wh{y}(\cdot)^*$ and $Y(\cdot)^* \wh{\rho}(\cdot)$ belong to $\sW_{-, 0}^{p\ts p}$, and thus $\om$ belongs to  $\sW_{-, 0}^{(p-m)\ts (p-m)}$. On the other hand, the function  $E^*Y(\cdot)^*Y(\cdot)E$ is hermitian on the imaginary axis, and hence the same is true for $\om$. But  for any positive integer $k$ we have
\[
\sW_{-,0}^{k\ts k}\cap (\sW_{-,0}^{k\ts k})^*=\{0\}.
\]
Thus $\om$ is identically zero, and thus $\tht(s)^*\tht(s)=E^*Y(s)^*Y(s)E=I_{p-m}$ for any $s\in i\BR$. Moreover, $\tht(\infty)^*\tht(\infty)=E^*E=I$.  This proves that $\tht$ is inner.
\epr

\begin{lem}\label{lem:GY3} Assume that  $T_G$ is right invertible, and   let $Y\in \sW_+^{p\ts p}$ be the function defined by \eqref{defY2}. Then $Y$ is  invertible in $\sW_+^{p\ts p}$.
\end{lem}
\bpr Fix $s\in i\BR$, and assume  $u\in\BC^p$ such that  $Y(s)u=0$. Then $G(s)Y(s)=D$ implies that $Du=0$.  By definition of $E$,  $u=Ev$ for some $v\in \BC^{p-m}$. Next use  $\tht(s)=Y(s)E$. It follows that $\tht(s)v=Y(s)Ev=Y(s)u=0$.   However, $\tht(s)$ is an isometry, by Lemma  \ref{lem:Theta}. So $v=0$, and  hence $u=0$. We see that $\det Y(s)\not =0$. Also $Y(\iy)=I_p$. We conclude that $T_Y$ is a Fredholm operator;  see \cite[Theorem XII.3.1]{GGK1}.

Next we prove that $\kr T_Y=\{0\}$. Take $h\in \kr T_Y$. Then $T_Yh=0$, and  hence  $Y(s)\wh{h}(s)=0$ for each $s\in i\BR$.  But $\det Y(s)\not =0$ for each $s\in i\BR$. Hence $\wh{h}=0$, and therefore $h=0$.

We want to prove that  $T_Y$ is invertible. Given the results of the preceding first two paragraphs it  suffices  to show that $\ind T_Y=0$. This will be done in the next step by an approximation argument, using the fact, from \cite{FKR4}, that we know the result is true for rational matrix functions.

Let $g$ be as in \eqref{GWiener}. Note  that $g$ is the limit in $L^1$ of a sequence $g_1, g_2, \ldots$ such that $G_n(s)=D+\wh{g_n}(s)$ is a stable rational matrix function; cf., Part (v) on page 229 of \cite{GGK1}. Since $T_G$ is right invertible, $T_{G_n}$ will also be right invertible for $n$ sufficiently large. In fact, $T_{G_n}T_{G_n}^*\to T_GT_G^*$ in operator norm. Put $y_n=T_{G_n}^*(T_{G_n}T_{G_n}^*)^{-1} g_n$. Then $y_n\to y$ in the $L^1$-norm. Put $Y_n(s)=I-\wh{y_n}(s)$. Then $T_{Y_n}\to T_Y$ in operator norm. For $n$ sufficiently large the operator $T_{Y_n}$ is invertible (see the paragraph preceding Theorem 1.2 in \cite{FKR4} and formula (2.17) in \cite{FKR4}). In particular, the Fredholm index  of  $T_{Y_n}$ is zero. But $\ind T_Y=\lim_{n\to \iy} \ind T_{Y_n}=0$. Thus $T_Y$ is invertible, and hence $Y$ is invertible in  $\sW_+^{p\ts p}$.\epr

\smallskip
\noindent\textbf{Proof of  Theorem \ref{mainthm2}.}   From  Lemma \ref{lem:GY3} we know that $Y\in \sW_+^{p\ts p}$ and that $Y$  is invertible in $\sW_+^{p\ts p}$. Thus we only have to prove the second part of the theorem.   Since $Y$ is invertible in $\sW_+^{p\ts p}$, the  $p\ts p$ matrix function given by \eqref{out44}  is well-defined and  belongs to $ \sW_+^{p\ts p}$. Furthermore, from  \eqref{idDE} we know that   the $p\ts p$ matrix $\begin{bmatrix}D^+ & E\end{bmatrix}$ is invertible. Hence the function   defined by the right hand side of  \eqref{out4} belongs to $\sW_+^{p\ts p}$ and is invertible in  $\sW_+^{p\ts p}$.    Using \eqref{propY1} and the identity \eqref{idDE} we see that
\begin{align*}
& \begin{bmatrix}
  G(s) \\
  E^* Y(s)^{-1}
\end{bmatrix}Y(s) \begin{bmatrix}D^+ & E\end{bmatrix}= \begin{bmatrix} G(s)Y(s) \\ E^* \end{bmatrix} \begin{bmatrix}D^+ & E\end{bmatrix}\\
&\hspace{1cm}=\begin{bmatrix} D \\ E^* \end{bmatrix} \begin{bmatrix}D^+ & E\end{bmatrix}=\begin{bmatrix} I_m & 0\\ 0& I_{p-m}\end{bmatrix}, \quad \Re s\geq 0.
\end{align*}
This proves the first identity  \eqref{out4}. The second identity is an immediate consequence of the definitions of $\Xi$ and $\tht$ in \eqref{defXi} and \eqref{defTheta}, respectively.\epr

\setcounter{equation}{0}
\section{Proof  of Theorem \ref{mainthm1}}\label{sec:prthm1}

We begin with  a lemma concerning  the functions $\Xi$ and  $\tht$.

\begin{lem}\label{lem:Y4}
Assume that  $T_G$ is right invertible, and let $\Xi$ and  $\tht$ be the functions defined by \eqref{defXi} and \eqref{defTheta}, respectively.  Then
\begin{equation}\label{eq:XiTht}
\kr T_G = T_\tht L^2_{p-m}(\BR_+)\quad\mbox{and}\quad \tht^*\Xi\in \wt{\sW}_{0, -}^{\,(p-m)\ts m}.
\end{equation}
\end{lem}

\bpr
We split the proof into two parts.

\smallskip
\noindent\textsc{Part 1.} In this part we prove the inclusion of \eqref{eq:XiTht}.
Take $s\in i\BR$. Note that in Proposition \ref{lem:Theta} it was shown for $s\in i\BR$ that
\[
Y(s)=I-G(s)^*  \wh{f}(s) + \wh{\rho}(s),
\]
where $f=(T_GT_G^*)^{-1}g$ and $\rho$ is defined by \eqref{defrho}; see \eqref{eqY3}.
From \eqref{defXi} and \eqref{defTheta}  we then see  that
\[
\tht(s)^*\Xi(s)=E^*Y(s)^*Y(s)D^+=E^*Y(s)^*\Big(I-G(s)^*  \wh{f}(s) + \wh{\rho}(s)\Big)D^+.
\]
Now use that $G(s)Y(s)E=DE=0$, and hence $E^*Y(s)^*G(s)^*=0$. The latter identity and the fact that $E^*D^+=0$ and $Y = I-\wh{y}$ imply  that
\begin{align}
\tht(s)^*\Xi(s)&=-E^*\wh{y}(s)^* D^+ +  E^*\wh{\rho}(s)D^+-E^*\wh{y}(s)^*\wh{\rho}(s)D^+\nn \\
&=-A(s)+B(s)-C(s).\label{defABC}
\end{align}
From item (b) in Proposition \ref{lem:Y4} we know that $y\in L_{p\ts p}^1(\BR_+)\cap  L_{p\ts p}^2(\BR_+)$, and thus  $\wh{y}\in \wt{\sW}_+^{\,p\ts p}$ and $\wh{y}(\iy)=0$, that is, $\wh{y}\in \wt{\sW}_{0,+}^{\,p\ts p}$. It follows that
\begin{equation}\label{propA}
A(\cdot):= E^*\wh{y}(\cdot)^* D^+ \in   \wt{\sW}_{0,-}^{\,(p-m)\ts p}.
\end{equation}
Recall that $\rho$ is given by \eqref{defrho} with $f=(T_GT_G^*)^{-1}g$.  Since the function $g$ belongs to  $L_{m\ts p}^1(\BR_+) \cap L_{m\ts p}^2(\BR_+)$, item (b) in Proposition \ref{prop:HGR} tells us that the same holds true for $f$.  It follows that $g^*\star f \in L_{p\ts p}^1(\BR)\cap  L_{p\ts p}^2(\BR)$. The latter implies that  $\rho\in L_{p\ts p}^1(\BR_-)\cap  L_{p\ts p}^2(\BR_-)$. We conclude that
\begin{equation}\label{propB}
B(\cdot):=E^*\wh{\rho}(\cdot)D^+\in   \wt{\sW}_{0,-}^{\,(p-m)\ts p}.
\end{equation}
Finally, note that $y^*(t)=y(-t)^*$ for $t \in \BR$ and $\big(\wh{y}(s)\big)^*=\wh{y^*}(s)$ for $s\in i\BR$.  Thus
\[
\wh{y}(s)^*\wh{\rho}(s)=\Big( \widehat{y^*\star \rho}\Big)(s) , \quad  s\in i\BR,
\]
and
\[
(y^*\star \rho)(t)= \int_{-\iy}^\iy y^*(t-\tau)\rho(\tau)\, d\tau=\int_{-\iy}^0 y^*(t-\tau)\rho(\tau)\, d\tau.
\]
Since both $y^*$ and $\rho$ belong to  $L_{p\ts p}^1(\BR_-)\cap  L_{p\ts p}^2(\BR_-)$,  it is well known (see, e.g., Section 2 in \cite{DG80}) that the same holds true for $y^* \star \rho$. But  then
\begin{equation}\label{propC}
C(\cdot):=E^*\wh{y}(\cdot)^*\wh{\rho}(\cdot)D^+\in  \wt{\sW}_{0,-}^{\,(p-m)\ts p}.
\end{equation}
From \eqref{propA},  \eqref{propB}, \eqref{propC} and  \eqref{defABC} it follows that  $\tht^*\Xi\in \wt{\sW}_{0, -}^{\,(p-m)\ts m}$.

\smallskip
\noindent\textsc{Part 2.} In this part we prove  the identity of \eqref{eq:XiTht}.
 Using \eqref{propY1} we see that
\[
G(s)\tht(s)=G(s)Y(s)E=DE=0, \quad s\in i\BR.
\]
This implies that $T_GT_\tht=0$, and hence  $\im T_\tht\subset \kr T_G$. To prove the reverse inclusion, take $h\in \kr T_G$. Thus $h\in L_p^2(\BR_+)$  and $T_Gh=0$.  It follows that $G(s)\wh{h}(s)=0$ for $\Re s >0$. Put $H(s)=\wh{h}(s)$. Then $H(\cdot)$ belongs to $H_m^2(i\BR)$. Next we apply Theorem \ref{mainthm2}. Using the identities  in \eqref{out4} we see that
\begin{align*}
   H(s)&= \begin{bmatrix}\Xi(s)&\tht(s)\end{bmatrix}\begin{bmatrix}
  G(s) \\ E^* Y(s)^{-1} \end{bmatrix}H(s)\nonumber  \\
 &= \begin{bmatrix}\Xi(s)&\tht(s)\end{bmatrix}\begin{bmatrix}
0 \\ E^* Y(s)^{-1} H(s)   \end{bmatrix}= \tht(s)E^* Y(s)^{-1}H(s).
\end{align*}
Hence  $\wh{h}(s)=\tht(s)\Psi(s)$, where $\Psi(s)= E^*Y(s)^{-1}\wh{h}(s)$. Since $h\in L_p^2(\BR_+)$  and $Y(\cdot)^{-1}$ is a matrix function with $H^\iy$ entries, we conclude that $\Psi\in H_{p-m}^2$, and hence $\Psi=\wh{u}$ for some $u\in  L_{p-m}^2(\BR_+)$.  The identity $\wh{h}(s)=\tht(s)\Psi(s)$ then yields $\wh{h}(s)=\tht(s)\wh{u}(s)$. This shows that   $h=T_\tht u$, and  hence $\kr T_G\subset  \im T_\tht$.
\epr

\medskip
\noindent\textbf{Proof of Theorem \ref{mainthm1}.} From Lemma \ref{lem:Theta} we know that  $\tht$ is an inner function in $\sW_+^{(p-m)\ts m}$. The   proof of the  other statements is split  into  three  parts.

\smallskip
\noindent\textsc{Part 1.} In this part we show that  the equation $G(s) X(s)  = I_m$, $\Re s > 0$,  has a solution $X\in \sW_+^{p\ts m}$  if and only if $T_G$ is right invertible. Furthermore, we show that in that case  the function $\Xi$  defined by \eqref{defXi} is a  particular solution.  From the one but last sentence  of the paragraph containing \eqref{Bez1} we know  that  it suffices to prove the ``if part'' only. Therefore, in what follows we assume that $T_G$ is right invertible. Since  $\Xi(s)=Y(s)D^+$ and $Y \in \sW_+^{p\ts p}$, we have $\Xi \in \sW_+^{p\ts m}$. Moreover, using the identity \eqref{lem:Y4} we have
\[
G(s)\Xi(s)=G(s)Y(s)D^+=DD^+=I_m.
\]
 Thus $\Xi$   is a particular solution.

\smallskip
\noindent\textsc{Part 2.}  This  second part deals with the description of all in solutions in $\sW_+^{p\ts m}$. Let $Z$ be an arbitrary function in $\sW_+^{p\ts m}$, and let $X\in \sW_+^{p\ts m}$ be defined by \eqref{sol8}.  Then
\[
G(s)X(s)=G(s)\Xi(s)+G(s)\tht(s)Z(s)=I_m+G(s)\tht(s)Z(s), \quad \Re s\geq 0.
\]
Recall that $G(s)\tht(s)=G(s)Y(s)E=DE=0$.  Thus $G(s)X(s)=I_m$, $\Re s \geq 0$, and thus $X$ is a solution.

To prove the converse implication, let $X\in \sW_+^{p\ts m}$  be a solution   of the equation $G(s)X (s)=I_m$. Put $H=X-\Xi$. Then $H\in \sW_+^{p\ts m}$  and $G(s)H (s)=0$.  Using  the identities in \eqref{out4}, we obtain
\begin{align*}
H(s)&=\begin{bmatrix}\Xi(s)&\tht(s)\end{bmatrix}\begin{bmatrix}
  G(s) \\ E^* Y(s)^{-1} \end{bmatrix}H(s)\nn  \\
 &= \begin{bmatrix}\Xi(s)&\tht(s)\end{bmatrix}\begin{bmatrix}
0 \\ E^* Y(s)^{-1}H(s)  \end{bmatrix}= \tht(s)E^* Y(s)^{-1}H(s).
\end{align*}
Thus $H(s)= \tht(s)Z(s)$, where $Z(s)=E^* Y(s)^{-1}H(s)$. Since $Y$ is invertible in $\sW_+^{p\ts p}$, the function  $Y(\cdot)^{-1}$ is in  $\sW_+^{p\ts p}$ . Together with the fact that $H\in \sW_+^{p\ts m}$, this yields  $Z\in \sW_+^{(p-m)\times m}$. It follows $X$ has the desired representation \eqref{sol8}.

\smallskip
\noindent\textsc{Part 3.} In this part we prove the identity \eqref{mincond}. Assume $Z\in \wt{\sW}^{\, (p-m)\ts m}$, and  let  $X$ be the function defined  by \eqref{sol8}. Fix $u\in \BC^m$. Then $Z(\cdot)u\in H_m^2(i\BR)$, and
\begin{equation}\label{normZ}
\tht(\cdot)Z(\cdot)u=M_\tht Z(\cdot)u\in H_{p-m}^2(i\BR).
\end{equation}
Here $M_\tht$ is the operator of multiplication by $\tht(\cdot)$ mapping  $H_m^2(i\BR)$ into $H_{p-m}^2(i\BR)$. Furthermore, since $M_\tht$ is an isometry,  we also see that
\begin{equation}
\label{normZ2}
\|Z(\cdot)u\|= \|\tht(\cdot)Z(\cdot)u\|.
\end{equation}
The fact that $y\in L^1_{p\ts p}(\BR_+)\cap L^2_{p\ts p}(\BR_+)$  implies that $Y\in \wt{\sW}_+^{\, p\ts p}$. But then $\Xi(s)=Y(s)D^+$ yields  $\Xi(\cdot)u\in H_m^2(i\BR)$. Using the identity  \eqref{sol8} we conclude that $\Xi(\cdot)u$ also belongs to $H_m^2(i\BR)$. It follows that all norms in \eqref{mincond} are well defined, and in order to prove the identity \eqref{mincond} it suffices that to show that in $H_m^2(i\BR)$ the function  $\tht(\cdot)Z(\cdot)u$  is orthogonal to the function  $\Xi(\cdot)v$ for any $v\in\BC^m$.  The latter fact follows from  the inclusion in the second part of \eqref{eq:XiTht}.  Indeed,  this inclusion  tells us that $M_\tht^*\Xi(\cdot) v=0$, and hence
\begin{align*}
\lg \Xi(\cdot) {v}, \tht(\cdot)Z(\cdot)u\rg_{H_m^2(i\BR)}
&=\lg \Xi(\cdot) {v}, M_\tht Z(\cdot)u\rg_{H_m^2(i\BR)}\\[.1cm]
&=\lg M_\tht^*\Xi(\cdot) {v}, Z(\cdot)u\rg_{H_m^2(i\BR)}=0.
\end{align*}
This completes the proof. \epr


\appendix
\section{The Lebesgue space  $L ^1(\BR)\cap L ^2(\BR)$}\label{sec:App}
\renewcommand{\theequation}{A\arabic{equation}}
\setcounter{equation}{0}
The material in this section is standard and  is presented for the convenience of the reader.   Throughout we deal with   the Lebesgue spaces  of complex-valued functions on the real line  $L ^1(\BR)$ and $L ^2(\BR)$, their vector-valued counterparts $L_m^1(\BR)$ and $L_m^2(\BR)$, and the intersection of the latter two spaces: $L_m^1(\BR)\cap L_m^2(\BR)$. The norms on these spaces  are given by
\begin{align*}
\|f\|_1    &=\int_{-\iy}^\iy |f(t)|\, dt \qquad \mbox{for $f\in  L^1(\BR)$,}  \\
\|f\|_2    &=  \Big(\int_{-\iy}^\iy |f(t)|^2\, dt\Big)^{1/2} \qquad \mbox{for $f\in  L^2(\BR)$},\\
\|f\|_1    &=\Big(\sum_{i=1}^m \|f_i\|_1^2\Big) ^{1/2} \qquad \mbox{for $f=(f_1, \ldots, f_m)^\top \in  L_m^1(\BR)$},  \\
\|f\|_2    &= \Big(\sum_{i=1}^m \|f_i\|_2^2\Big)^{1/2} \qquad \mbox{for $f=(f_1, \ldots, f_m)^\top \in  L_m^2(\BR)$},\\
\|f\|_0    &= \max\{\|f\|_1, \|f\|_2\} \qquad \mbox{for $f\in  L_m^1(\BR)\cap L_m^2(\BR)$}.
\end{align*}

Let $k\in L_{m\ts p}^1(\BR)$. Thus $k$ is an $m\ts p$ matrix function of which the $(i,j)$-th entry $k_{ij}\in L^1(\BR)$.  For each $\va\in L_p^1(\BR) $ and $\psi \in L_p^2(\BR)$  the convolution products $k\star \va$ and $k\star \psi$,  see \eqref{conprod},  are  well defined, $k\star \va$ belongs to $L_m^1(\BR)$ and $k\star \psi$ belongs to $L_m^2(\BR)$. In particular, if  $f\in L_p^1(\BR)\cap L_p^2(\BR)$, then $k\star f$ belongs to $L_m^1(\BR)\cap L_m^2(\BR)$.  It follows that for a given $k\in L_{m\ts p}^1(\BR)$ the convolution product induces linear maps from the space  $L_p^1(\BR)$ into $L_m^1(\BR)$,   from the space  $L_p^2(\BR)$ into $L_m^2(\BR)$, and  from  the space  $L_p^1(\BR)\cap  L^2_p (\BR)$ into $L_m^1(\BR)\cap L_m^2(\BR)$. The resulting operators will be denoted by  $K_1$, $K_2$ and $K_0$, respectively.  The proof of the following lemma is standard (see, e.g.,  page 216 in [6]) and therefore it  is  omitted.

\begin{lem}\label{lemK} The  operators $K_1$, $K_2$ and $K_0$ are bounded linear operators, and
\begin{equation}
\label{bdnorm}
\|K_\nu\|\leq \kappa\hspace{.2cm}(\nu=1, 2, 0), \hspace{.2cm} \mbox{where}\hspace{.2cm}  \kappa=\Big(\sum_{i=1}^m \sum_{j=1}^p\|k_{ij}\|_1^2\Big)^{1/2}.
\end{equation}
\end{lem}

With $k\in L_{m\ts p}^1(\BR_+)$ we also associate  the Wiener-Hopf operator $W$ and the Hankel operator $H$ defined by
\begin{align*}
(Wf)(t)    &=\int_{0}^\iy  k(t-\t)f(\t)\,d\t, \quad 0\leq t<\iy,\\
(Hf)(t)&=\int_{0}^\iy  k(t+\t)f(\t)\,d\t, \quad 0\leq t<\iy.
\end{align*}
Using the classical relation  between the convolution operator defined by $k$ and
the operators $W$ and $H$ (see, e.g.,  Section XII.2 in \cite{GGK1})  it is easy to see that $W$ and $H$  map   the  space $L_p^1(\BR_+)$ into $L_m^1(\BR_+)$,   the space  $L_p^2(\BR_+)$ into $L_m^2(\BR_+)$, and  the space  $L_p^1(\BR_+)\cap L_p^2(\BR_+)$ into $L_m^1(\BR_+)\cap L_m^2(\BR_+)$.  We denote the resulting operators by $W_1$,  $W_2$, $W_0$, and  $H_1$,  $H_2$, $H_0$, respectively. Lemma \ref{lemK} shows that these operators are bounded and
\begin{equation}\label{bdnormWH}
\|W_\nu\| \leq \k\hspace{.1cm} \mbox{and}\hspace{.1cm} \|H_\nu\| \leq \k\hspace{.2cm}(\nu=1, 2, 0), \hspace{.1cm}\mbox{where} \hspace{.1cm}\k=\Big(\sum_{i=1}^m \sum_{j=1}^p\|k_{ij}\|_1^2\Big)^{1/2}.
\end{equation}
Furthermore, using the line of reasoning in Lemma XX.2.4 in \cite{GGK1},  we have the following corollary.

\begin{cor}\label{corH} The Hankel operators $H_1$, $H_2$, and $H_0$ are the limit in operator norm of  finite rank operators, and hence compact.
\end{cor}

Next we present an auxiliary result  that is  used in  the  proof of Proposition~\ref{prop:HGR}.  Put
\begin{equation}
\label{RWiener}
R(s)=D_R+ \int_{-\iy} ^\iy e^{-st} r(t)\, dt \quad\mbox{where}\quad  r\in L_ {m\ts m}^1(\BR).
\end{equation}
By $T_R$ we denote the Wiener-Hopf operator on $L_m^2(\BR_+)$ defined by $R$, that is,
\begin{equation}
\label{TWiener}
 (T_Rf)(t)= D_Rf(t)+\int_{0}^\iy  r(t-\t)f(\t)\,d\t, \quad 0\leq t<\iy.
 \end{equation}
As we know from the first paragraph  of this section,  the fact that   $r\in L_{m\ts m}^1(\BR)$ implies that  $T_R$ maps  $L_m^1(\BR_+)\cap L_m^2(\BR_+)$ into itself.

\begin{lem}\label{lemR}
If  $T_R$ is invertible as an operator on $L_m^2(\BR_+)$, then $T_R^{-1}$  maps  the space $L_m^1(\BR_+)\cap L_m^2(\BR_+)$ in a one-to-one way onto itself.
\end{lem}
\bpr Since $T_R$ is invertible, $R$ admits a canonical factorization (see Section XXX.10 in \cite{GGK2}), and hence we can write $T_R^{-1}=LU$, where $L$ and $U$ are Wiener-Hopf operators on $L_m^2(\BR_+)$,
\begin{align}
(Lf)(t)   &=   D_L f(t)+\int_0^t \ell(t-\t)f(t)\, dt,\quad  0\leq t<\iy,  \label{eqL}\\
(Uf)(t)  &=  D_U f(t)+\int_t^\iy  u(t-\t)f(t)\, dt, \quad  0\leq t<\iy.  \label{eqU}
\end{align}
Here   $\ell$ and $u$ both belong to $L_{m\ts m}^1(\BR)$, with support of $\ell$ in $ \BR_+$ and  support of $u$ in  $\BR_-$. The fact that both $\ell$ and $u$ belong to $L_{m\ts m}^1(\BR)$ implies that both $L$ and $U$ map  $L_m^1(\BR_+)\cap L_m^2(\BR_+)$ into itself. Hence $T_R^{-1}$ has the same property.
Since $T_R^{-1}$ is one-to-one on  $L_m^2(\BR_+)$, it is also  one-to-one  on $L_m^1(\BR_+)\cap L_m^2(\BR_+)$. For $f\in L_m^1(\BR_+)\cap L_m^2(\BR_+)$, we have $g=T_R f\in L_m^1(\BR_+)\cap L_m^2(\BR_+)$ and $f= T_R^{-1}T_R f=T_R^{-1}g$. This shows that $T_R^{-1}$ maps $L_m^1(\BR_+)\cap L_m^2(\BR_+)$ onto $L_m^1(\BR_+)\cap L_m^2(\BR_+)$. \epr

\medskip

\begin{lem}\label{lem4}
Let  $k\in L_{p\ts m}^1(\BR_+)$ and $\tilde{k}\in L_{m\ts p}^1(\BR_+)$, and let $H$ and $\tilde{H}$ be the corresponding Hankel operators acting from $L_m^2(\BR_+)$ into $L_p^2(\BR_+)$ and from $L_p^2(\BR_+)$ into $L_m^2(\BR_+)$, respectively.  Let $Q$ be any operator on $L_m^2(\BR_+)$  mapping $L_p^1(\BR_+)\cap L_m^2(\BR_+)$ into itself, and assume that the restricted  operator $Q_0$ acting on $L_m^1(\BR_+)\cap L_m^2(\BR_+)$  is bounded. If the operator $I-\tilde{H}QH$ is   invertible  on $L_m^2(\BR_+)$, then  $I-\tilde{H}QH$ maps  the space $L_m^1(\BR_+)\cap L_m^2(\BR_+)$ in a one-to-one way onto itself.
\end{lem}
\bpr  We know that  $H$  maps $L_p^1(\BR_+)\cap L_2^p(\BR_+)$ into  $L_m^1(\BR_+)\cap L_m^2(\BR_+)$. Furthermore the same holds true for $\tilde{H}$ with the role of $p$ and $m$ interchanged.  Hence our hypothesis on $Q$ implies that $I-\tilde{H}QH$ maps  the space $L_m^1(\BR_+)\cap L_m^2(\BR_+)$ into itself. Let $M_0$ be the corresponding restricted operator.  We have to prove that $M_0$ is invertible.  Note that Corollary \ref{corH} implies that $M_0$ is equal to the identity operator minus a compact operator, and hence $M_0$ is a Fredholm operator of index zero. Therefore,  in order to prove that $M_0$ is invertible, it suffices to show that $\kr M_0$  consists of the zero element only. Assume not. Then there exists a non-zero $f$ in  $L_m^1(\BR_+)\cap L_m^2(\BR_+)$ such that  $M_0f=0$. The fact that $f$ belongs to $L_m^1(\BR_+)\cap L_m^2(\BR_+)\subset L_m^2(\BR_+)$ shows that $0=M_0f=(I- \tilde{H}QH)f$. But $I- \tilde{H}QH$ is assumed to be invertible. Hence $f$ must be zero.
Thus $M_0$ is invertible. \epr

\section*{Acknowledgement}
This work is based on the research supported in part by the National Research Foundation. Any opinion, finding and conclusion or recommendation expressed in this material is that of the authors and the NRF does not accept any liability in this regard.

\end{document}